\documentclass[10pt]{article}
\usepackage{amssymb}
\usepackage{latexsym}
\newtheorem{thm}{Theorem}[section]
\newtheorem{prop}[thm]{Proposition} 
 
 \newtheorem{rmk}[thm]{Remark}

\newcommand {\pf}{\noindent{\bf Proof.}\ }
\newcommand{\complex}{{\mathbb C}}

\newcommand{\reals}{{\mathbb R}}

\newcommand{\integers}{{\mathbb Z}}

\newcommand{\calg}{{\cal G}}

\newcommand{\calm}{{\cal M}}
\newcommand{\calo}{{\cal O}}

\newcommand{\calu}{{\cal U}}

\newcommand{\qed}{\begin{flushright} $\Box$\ \ \ \ \ \end{flushright}}

\newcommand{\arrows}{\,\lower1pt\hbox{$\longrightarrow$}\hskip-.24in\raise2pt
             \hbox{$\longrightarrow$}\,}

\newcommand{\Lag}{\mathrm{Lag}}
\newcommand{\topp}{\mathrm{top}}
\begin{document}

\title{{\bf The Maslov gerbe}}
\author
{Alan
Weinstein\thanks{Research partially supported by NSF Grant
DMS-0204100
\newline \mbox{~~~~}MSC2000 Subject Classification Number: 53D12
(Primary).
\newline \mbox{~~~~}Keywords: Maslov class, gerbe, lagrangian subspaces}\\
Department of Mathematics\\ University of California\\ Berkeley, CA
94720 USA\\ {\small(alanw@math.berkeley.edu)}}
\date{December 13, 2003}
\maketitle
\begin{abstract}
Let $\Lag(E)$ be the grassmannian of lagrangian subspaces of a complex
symplectic vector space $E$.  We construct a Maslov class which
generates the second integral cohomology of $\Lag(E)$, and we show
that its mod 2 reduction is the characteristic class of a flat gerbe
with structure group $\integers_2$.  We explain the relation of this
gerbe to the well-known flat Maslov line bundle with structure group
$\integers_4$ over the real lagrangian grassmannian, whose
characteristic class is the mod 4 reduction of the real Maslov class.
\end{abstract}
\section{Introduction}
\label{sec-intro}

In real symplectic geometry, the {\bf Maslov class} is a canonical
class $\mu(E)$ which generates the first integral cohomology of
 the lagrangian
grassmannian $\Lag(E)$ of a symplectic vector space $E$.  Given any
$L_0\in\Lag(E)$, the homology class dual to $\mu(E)$ is represented by
the {\bf Maslov cycle} of $E$ with respect to $L_0$; the latter is the
cooriented, codimension 1, real algebraic variety $\Sigma_{L_0}$ (with
singularities of codimension 3 in $\Lag(E)$) consisting of the
lagrangian subspaces in $E$ which are not transverse to $E$.  The mod
4 reduction of $\mu(E)$ has the further interpretation as the
characteristic class (or holonomy) of the {\bf Maslov line bundle}, a
flat $U(1)$ bundle over $\Lag(E)$.  Although the full Maslov class is
essential in symplectic topology, only the mod 4 reduction matters in
microlocal analysis, where the Maslov class was first discovered.
(See \cite{ma:theory}.  The topological description as a cohomology
class was first given by Arnol'd \cite{ar:characteristic}.)

When $E$ is the sum $V \oplus \overline V$ of a symplectic vector
space $V$ and the same vector space with the opposite
symplectic structure, $\Lag(E)$ contains the graphs of the symplectic
automorphisms of $V$, so the symplectic group $Sp(V)$ embeds in
$\Lag(E)$.  The Maslov class on $\Lag(E)$ pulls back to twice the
generator of $H^1(Sp(V),\integers)$, so a 
($4$-valued) parallel section of the Maslov line bundle
of $E$ restricts to $Sp(V)$ as a double covering which can be
identified with the metaplectic group $Mp(V)$.

This linear algebra and topology is applied to microlocal analysis as
follows \cite{ho:fourier}\cite{ma:theory}.  Over a cotangent bundle
$T^*X$, the tangent spaces to the fibres form a lagrangian subbundle
$L_0$ of the symplectic vector bundle $T(T^*X)$, i.e. a section of the
bundle $\Lag(T(T^*X))$ of lagrangian subspaces of the tangent spaces
to $T^*X$.  The Maslov cycles attached to the individual tangent
spaces then combine to form a cycle in $\Lag(T(T^*X))$.  This cycle is
dual to a canonical Maslov class in $H^1(\Lag(T(T^*X)),\integers)$,
whose mod 4 reduction is again the characteristic class of a flat
$U(1)$ bundle, the Maslov bundle $\calm_X\to \Lag(T(T^*X))$.

For each lagrangian immersion $j:\Lambda\to T^*X$, there is a
natural ``Gauss map'' $G_j:\Lambda \to \Lag(T(T^*X))$, and the
``Maslov canonical operator'' (or cousins going by various names)
produces (sometimes distributional, or depending on a ``quantization''
parameter) sections of the bundle $\sqrt{|\bigwedge^{\topp} T^*X|}$ of
half densities on $X$ from sections of $\sqrt{|\bigwedge^{\topp}
  T^*\Lambda|}\otimes G_j^*(\calm_X)$.  In other words, the 
Maslov bundle is a kind of ``transition object'' from the half
densities on lagrangian submanifolds of
$T^*X$ to those on $X$ itself.  

What becomes of all this when the real numbers are replaced by the
complex numbers?  In this paper, we will answer this question at the
level of linear algebra and topology, leaving applications to
microlocal analysis for the future. 

$\Lag(E)$ is now a complex manifold, and the Maslov cycle associated
to $L_0$
is a complex subvariety of complex codimension 1, with singularities of
complex codimension 3.  Repeating in the complex setting the proof of
Theorem 3.4.9 in 
\cite{du:fourier}, one can show 
that, as in the real case, the complement of the
Maslov cycle is contractible, and (since a complex subvariety is
coorientable) that this cycle is therefore dual to a generator of
$H^2(\Lag(E),\integers)$ which is independent of the choice of $L_0$.  We
call this generator the {\bf complex Maslov class}.  The main point of
this paper is that the mod 2 reduction of the complex Maslov class
is the characteristic class of a flat $U(1)$ {\em gerbe} over
$\Lag(E)$ which we naturally call the {\bf Maslov gerbe}.  Objects of
this gerbe are local square roots of the unitarized determinant 
bundle associated to the 
tautological vector bundle over $\Lag(E)$.  The absence of global
objects of this gerbe is well known and corresponds to the
nonexistence of various other objects, including half-forms on general
complex manifolds and a connected double covering of
$Sp(2n,\complex)$.  The latter impossibility plays a central role in
recent work of Omori et al \cite{om-ma-mi-yo:strange}
 on star-exponentials of complex
quadratic polynomials, which was in fact the stimulus for the present
work.

Much of the work described here was done while I was a fellow of the
Japan Society for the Promotion of Science at
Keio University.  I thank the JSPS and my host, Yoshiaki Maeda, for
the opportunity to make this visit.  The work was continued while I was
visiting the Institut Math\'ematique de Jussieu.  I would like to
thank Joseph Oesterl\'e and Harold Rosenberg for their hospitality,
and Pierre Schapira for encouragement to write this note.  Throughout
this work, I have been especially stimulated by discussions with Maeda
and with Pedro Rios.  I also thank Ralph
Cohen, Tara Holm, and Sergey Lysenko for their helpful comments.

\section{Square roots of line bundles as $\integers_2$ gerbes}
\label{sec-squareroots}
Let $\lambda$ be a complex 
line bundle over a manifold $M$.  Following Brylinski \cite{br:loop},
we may define a $\integers_2$ gerbe $\sqrt{\lambda}$ over $M$ whose
nontriviality is the obstruction to the existence of a (tensor) square
root of $\lambda$.  Namely, we apply Proposition 5.2.3 of
\cite{br:loop} to the exact sequence of groups 
$$0\to \integers_2\to GL(1) \stackrel{s}{\to} GL(1)\to 0,$$
 where $s$ is the squaring
homomorphism.  (We will identify $\integers_k$ with
the multiplicative group of $k$'th roots of unity.)  In the case of a
hermitian line bundle, everything is essentially the same, with
the sequence above replaced by 
$$0\to \integers_2\to U(1) \stackrel{s}{\to} U(1)\to 0.$$

According to the cited proposition, there is associated to the 
line 
bundle $\lambda$ a $\integers_2$ gerbe which represents the
obstruction to lifting the bundle through $s$.  A
construction of Giraud (Theorem 5.2.8 in \cite{br:loop}) associates to
this gerbe a characteristic class $\gamma(\sqrt{\lambda})\in
H^2(M,\integers_2)$.  Finally, Theorem 5.2.9 in \cite{br:loop}
asserts that this characteristic class is obtained from the
characteristic class of $\lambda$ in sheaf cohomology by the
connecting homomorphism of one of the exact sequences above.  Comparing these
with the exponential sequences usually associated with line bundles, 
$$0\to\integers\to \complex\to GL(1)\to 0,$$
or 
$$0\to\integers\to \reals\to U(1)\to 0$$
in the hermitian case,
we conclude that the Giraud
characteristic class is simply the mod 2 reduction of the Chern
class $c(\lambda)$.

These constructions may be described in more concrete terms.
For each open $\calu\subseteq M$,
$\sqrt{\lambda}(\calu)$ is the groupoid whose objects are pairs
$(\tau,\iota)$ consisting of a (hermitian in the $U(1)$ case)
 line bundle $\tau$ and 
an isomorphism $\iota$ from the tensor square $\tau^2$ to the restriction
$\lambda |_\calu$.  A morphism from $(\tau,\iota)$ to $(\tau',\iota)$
is a bundle isomorphism $\sigma:\tau\to\tau'$ such that $\iota'
\sigma^2 \iota^{-1}$ is the identity automorphism of $\lambda
|_\calu$, where $\sigma^2$ is the tensor square of $\sigma$.  It is
easy to see that any two objects in $\sqrt{\lambda}(\calu)$ are
isomorphic and that the automorphism group of $(\tau,\iota)$ may be
identified with the continuous (hence locally constant) functions on
$\calu$ with values in $\integers_2$.  Thus, $\sqrt{\lambda}$ is a
gerbe with band $\integers_2$.

To determine the Giraud characteristic class of a gerbe $\calg$ with
discrete structure group, we 
choose a good covering of $M$ by open subsets $\calu_i$ which are the
domains of objects $\calo_i$.  On $\calu_i \cap \calu_j$, we
choose an isomorphism $\sigma_{ij}$ from the restriction of
$\calo_j$ to that of $\calo_i$.  On a triple
intersection $\calu_i\cap\calu_j\cap \calu_k$, the composition
$\gamma_{ijk}=\sigma_{ij}\sigma_{jk}\sigma_{ki}$ is a constant
function with values in the structure group; these functions define a \v{C}ech
2-cocycle which is, up to coboundaries, independent of choices.  We
denote the resulting cohomology class by $\gamma(\calg)$.

To relate the Giraud class of $\sqrt{\lambda}$ to the Chern class
of $\lambda$, we determine the
Chern class by choosing a good covering of $M$ by open subsets
$\calu_i$ and trivializations $\epsilon_i:\lambda|_{\calu_i} \to
\calu_i \times \complex$.  On $\calu_i \cap \calu_j$, the composition
$\epsilon_i\epsilon_j^{-1}$ is represented by a $GL(1)$-valued function
$r_{ij}$.  We choose complex-valued (real-valued in the hermitian
case) 
functions $\theta_{ij}$ such that
$r_{ij}=e^{2\pi i \theta_{ij}}.$ On $\calu_i\cap\calu_j\cap \calu_k$,
$r_{ij}r_{jk}r_{ki}=1,$ so $c_{ijk}=
\theta_{ij}+\theta_{jk}+\theta_{ki}$ is a constant $\integers$-valued
function.  The integers $c_{ijk}$ form a cocycle which represents the
Chern class $c(\lambda)\in H^2(M,\integers)$.

We now compute the Giraud class by using the
same covering and letting $\tau_i=\calu_i
\times \complex$, with $\iota_i((m,z)\otimes(m,z))
=\epsilon_i^{-1}(m,z^2).$  Choosing
$\sigma_{ij}(m,z)=(m,e^{\pi i \theta_{ij}(m)}z)$ gives
$$\gamma_{ijk}=\sigma_{ij}\sigma_{jk}\sigma_{ki}=e^{\pi i
(\theta_{ij}+\theta_{jk}+\theta_{ki})}=e^{\pi i
c_{ijk}}=(-1)^{c_{ijk}}.$$ 

Thus we have proved:

\begin{prop}
\label{prop-sqrt}
If $\lambda$ is a complex line bundle over $M$, then the Giraud
class $\gamma(\sqrt{\lambda})$ is the mod 2 reduction of the Chern 
class $c(\lambda)$.
\end{prop}

\section{Reduction modulo 2 of  the complex Maslov class}
\label{sec-maslovcomplex}

In this section, we will identify a line bundle $\lambda$ over
$\Lag(E)$ whose Chern class is the complex Maslov
class when $E$ is a complex symplectic vector space.

First, we note that, just as we can form tensor powers of line
bundles, 
we can also pass in a canonical way from a line bundle
$\lambda$ to a unitary line bundle $\arg(\lambda)$ which we will call the {\bf
  unitarization} of $\lambda$.  It is just the bundle 
associated to the principal bundle of frames of $\lambda$ via the homomorphism
$\arg:z\mapsto z/|z|$ from $GL(1,\complex)$ to $U(1)$.  
It is the tensor quotient of $\lambda$ by the
line bundle $|\lambda|$ obtained via the absolute value homomorphism
from $GL(1,\complex)$ to the positive real numbers $\reals_+$.  (When $\lambda$
is a real line bundle, we work with its complexification.)
Since
$\reals_+$ is
contractible, $|\lambda|$ is always trivial, and hence the characteristic
(Chern or Stiefel-Whitney) 
class of $\arg(\lambda)$ is the same as that of $\lambda$.  

Now let $\eta$ be the tautological vector bundle over $\Lag(E)$, i.e. the
bundle whose fibre over $L$ is $L$ itself.  
Fix a lagrangian subspace $L_0$ and 
a nonzero element $\phi$ of $\bigwedge^\topp (E/L_0)^*$.  The pullbacks of
$\phi$ to the elements of $\Lag(E)$ by the projection $E\to E/L_0$ define
a holomorphic  section of the line bundle
 $\bigwedge^{\topp}\eta^*$ which vanishes
precisely along the Maslov cycle $\Sigma_{L_0}$; it is transverse to
the zero section along the regular part of $\Sigma_{L_0}$, where
$L\cap L_0$ is 1-dimensional.  (Thus, the Maslov cycle is a divisor
associated to this line bundle.)  It follows immediately that the complex
Maslov class is the Chern class of the line bundle
$\bigwedge^{\topp}\eta^*$ and hence that of
$\arg(\bigwedge^{\topp}\eta^*)$ as well.  Applying Proposition \ref{prop-sqrt}, we
arrive at

\begin{prop}
\label{prop-gerbe}
Let $E$ be a (finite-dimensional) complex symplectic vector space.
The mod 2 reduction of the complex Maslov class is the Giraud
characteristic class of the $\integers_2$ gerbes
$\sqrt{\bigwedge^{\topp}\eta^*}$ and $\sqrt{\arg(\bigwedge^{\topp}\eta^*)}$, 
where $\eta$ is the tautological bundle
over the lagrangian grassmannian $\Lag(E)$, and $\arg(~)$ is the unitarization.
\end{prop}

\begin{rmk}
{\em We can say a bit more about the complex Maslov class.  The conormal
bundle to the Maslov cycle along the smooth locus is the complex line
bundle whose fibre over each element $L$ consists of the quadratic
forms on the line $L\cap L_0$.  Being the square of a line bundle,
this conormal bundle has a Chern class which vanishes mod 2.  But this
Chern class is dual to the self-intersection of the Maslov cycle with
itself, and so is equal to the square of the Maslov class.  Were it
not for the singularities of the Maslov cycle, we could 
conclude immediately that the cup square (in degree 4)
of the complex Maslov class is always an even class (i.e. zero modulo
2) or, equivalently, that the square of the (mod 2) Giraud class is zero.  
In fact, Ralph Cohen (private communication)
has verified this conclusion by
representing $\Lag(\complex ^{2n})$ as the homogeneous space
$Sp(n)/U(n)$.  The vanishing of the square of the complex Maslov class
 is also consistent with results in
Section 6 of \cite{go-ho:real} to the effect that the even $\integers_2$
cohomology rings of certain complex manifolds are isomorphic to the
$\integers_2$
cohomology rings of their real forms, via an isomorphism which halves
degrees.}
\end{rmk}

\section{Reduction modulo 4 of the real Maslov class}
\label{sec-maslovreal}

In preparation for our comparison of the real and complex Maslov classes, we give a geometric description of the Maslov bundle in the
real case as a square root of the flat $\integers_2$ bundle 
$\arg(\bigwedge^{\topp}\eta^*)$.

Since the trivial bundle $|\bigwedge^{\topp}\eta^*|$
has a natural ``positive'' square root, the bundle
$\sqrt{|\bigwedge^{\topp}\eta^*|}$ of half-densities, this is equivalent
to studying square roots of $\bigwedge^{\topp}\eta^*$, which are
generally known as bundles of half-forms.  (See, for instance,
\cite{gu-st:geometric}).

We look first at the case where $E$ is 2-dimensional, so that the
lagrangian grassmannian is simply a projective line, and
$\bigwedge^{\topp}\eta^*$ is just $\eta^*$ itself.  In terms of canonical
coordinates $(q,p)$ on $E$, we have two coordinate systems on
$\Lag(E)$.  The first assigns to each line of the form $p=aq$ its
slope $a$; the second assigns to each line of the form $q=bp$ its
inverse slope $b$.  The range of each coordinate system is an entire
(real or complex) line, and the transition map, defined where $a$ and
$b$ are nonzero, is $b=1/a$.  Bases for $\eta^*$ on the two
neighborhoods are given by $dq$ and $dp$, with the transition relation
$dp=a dq.$ Bases for the bundle $\sqrt{|\eta^*|}$ of half-densities
are then given by $\sqrt{|dq|}$ and $\sqrt{|dp|},$ with the transition
relation $\sqrt{|dp|}= \sqrt{|a|}\sqrt{|dq|}.$

Constructing a square root of $\arg(\eta^*)$ is equivalent to choosing a
square root  $\sqrt{\arg(a)}$ on the set where $a\neq 0$.
This can be done in two inequivalent ways, in the real case. 
For $a>0$, we may without loss of generality take the positive square root $1$, while for
$a<0$ we may
take either $i$ or ${-i}$.  Once we have made this
choice, we take bases for $\sqrt{\arg(\eta^*)}$ which we may call 
$\sqrt{\arg(dq)}$ and $\sqrt{\arg(dp)}$, with the transition relation
$\sqrt{\arg(dp)} = \sqrt{\arg (a)}\sqrt{\arg(dq)}$.

Now let us follow a parallel section of $\sqrt{\arg(\eta^*)}$ around a
loop in $\Lag(\reals^2)$ which starts at the line $a=1$ and rotates the line
in the counterclockwise direction in the $(q,p)$ plane.  As $a$
increases to infinity and the line becomes vertical, we may take
$\sqrt{\arg(dq)}$ as this section.   We write this section as 
$\sqrt{\arg(a)}^{-1}\sqrt{\arg(dp)}=\sqrt{\arg(dp)}$ and then follow it as the line crosses the $p$-axis into the region of negative $a$.  To cross the horizontal axis and return to the line $a=1$, we first return to the $q$ representation by writing $\sqrt{\arg(dp)}= \sqrt{\arg (a)}\sqrt{\arg(dq)}=\pm i \sqrt{\arg(dq)}$.  We have thus returned to $\pm i$ times our original section, and the holonomy around our loop is multiplication by $e^{\pm i \pi/2}$.  The sign depends on our choice of square root, which can be chosen so that this holonomy is given by the mod 4 reduction of the Maslov class.  We will not concern ourselves here with the choice of sign, which depends on the sign convention used in the definition of the Maslov class, via the choice of coorientation of the Maslov cycle.  

When the dimension of $E$ is greater than 2, we can write $E$ as a direct sum of a symplectic plane and another symplectic summand, and take a loop in $\Lag(E)$ consisting of the direct sum of a line in the symplectic plane, moving as in the paragraph above, with a fixed lagrangian subspace in the second summand.  This reduces the general case to that where $E$ is 2-dimensional, and we arrive at the following result.

\begin{prop}
\label{prop-flat}
Let $E$ be a (finite-dimensional) real symplectic vector space.  the mod 4 reduction of the real Maslov class is the characteristic
class of a (flat) $\integers_4$ bundle whose tensor square is  $\arg(\bigwedge^\topp \eta^*)$, where  $\eta$ is the tautological bundle over $\Lag(E)$.
\end{prop}

\section{Relation between the real and complex Maslov classes and
  their reductions} 
\label{sec-realcomplex}

We will begin with the following general result.

\begin{prop}
\label{prop-equator}
Let $\calg$ be a flat $U(1)$ gerbe over the sphere $\complex P^1$.  Let $\calo_+$ and $\calo_-$ be objects of $\calg$ defined on neighborhoods $\calu_+$ and $\calu_-$ of the northern and southern hemispheres respectively.  Then the value of the
 Giraud class of $\calg$ on the fundamental cycle of $\complex P^1$ is equal to the holonomy of the flat $U(1)$ bundle $\calo_+ \otimes \calo_-^{-1}$
 around the equator $\reals P^1$.
\end{prop}

\pf
We use the concrete construction of the Giraud class introduced in Section \ref{sec-squareroots}.
$\complex P^1$ may be triangulated as the boundary of a 3-simplex by four triangles, one of which exactly fills the southern hemisphere, while the other three have a common vertex at the north pole and are numbered in counterclockwise order when viewed from above.  We enlarge these triangles to open subsets $\calu_0 \ldots 
\calu_3$ which form a nice covering, and where $\calu_0=\calu_-$, while the other three subsets are contained in $\calu_+$.  We may then choose the object $\calo_0$ to be $\calo_i$, and the other three to be the restrictions of $\calo_+$.

When $i$ and $j$ are between $1$ and $3$, we choose the isomorphism
$\sigma_{ij}$ over $\calu_i \cap \calu_j$ to be the one induced from
the identifications with $\calo_+$.  We choose the isomorphisms
$\sigma_{0i}$ to be parallel sections of the flat bundles
$\calo_0\otimes \calo_i^{-1}$, all of which are identified with
$\calo_- \otimes \calo_+^{-1}$.  The value of the Giraud cocycle on
the fundamental class of the sphere is the product in $U(1)$ of the
four compositions (constant sections of trivial bundles)
$\gamma_{ijk}=\sigma_{ij}\sigma_{jk}\sigma_{ki}$ where $(ijk)$ takes
the values $(123)$, $(032)$, $(013)$, and $(021)$.  Now
$\sigma_{123}=1$ by construction, and, supressing the natural
isomorphisms, we find the remaining three values to be
$\sigma_{03}\sigma_{02}^{-1}$, $\sigma_{01}\sigma_{03}^{-1}$, and
$\sigma_{02}\sigma_{01}^{-1}$.  But the product of these three
``jumps" is exactly the inverse of the holonomy of the flat bundle
$\calo_- \otimes \calo_+^{-1}$ around the equator $\reals P^1$
oriented as the boundary of the northern hemisphere, i.e. with the
usual orientation on $\reals$ when the point at infinity is ignored.
Thus, it is equal to the holonomy of the inverse bundle $\calo_+
\otimes \calo_-^{-1}$.  
\qed

We apply the proposition above to the gerbe
$\sqrt{\bigwedge^{\topp}\eta^*}$ over $\Lag(E)$, where $E$ is complex.
The general case can be reduced to that where $E$ is 2-dimensional, so
we may identify $\Lag(E)$ with $\complex P^1$ as we have done before.
Objects $\calo_{\pm}$ on the northern and southern hemispheres must
restrict to the two different square roots of
$\bigwedge^{\topp}\eta^*$ on the equator.  (If they were the same,
there would be a global object.)  Therefore, the tensor quotient
$\calo_+ \otimes \calo_-^{-1}$ is the tensor square of one of these
square roots, i.e. of the Maslov line bundle or its inverse.  It is
for this reason (we might say that it is because the sphere is made of
{\em two} hemispheres), that the structure group $\integers_2 \subset
U(1)$ of the Maslov gerbe on the complex Lagrangian grassmannian is
the square of the structure group $\integers_4 \subset U(1)$ of the
Maslov line bundle on the real lagrangian grassmannian.

\begin{rmk}
{\em There is a similar phenomenon in ``one degree lower.''
In \cite{ho:fourier}, the holonomy in $\integers_4$ of
the Maslov line bundle is effectively expressed as the product of an
even number of jumps in $\integers_8$ (see the remark on page 163).}
\end{rmk}

\end{document}